\documentclass[12pt]{amsart}
\usepackage{amsthm,amsfonts,amsmath,amssymb,amsaddr,geometry,eufrak,bigints,eucal}
\usepackage{graphicx}
 \usepackage{hyperref}
\usepackage{enumitem}
\usepackage{cleveref}
\usepackage{orcidlink}
\usepackage{titlecaps}
\usepackage{xcolor} 

\hypersetup{
    colorlinks=true,       
    linkcolor=blue,        
    citecolor=red!100!black, 
    filecolor=magenta,     
    urlcolor=cyan          
}

\newtheorem{theorem}{Theorem}[section]
\newtheorem{lemma}[theorem]{Lemma}
\newtheorem{proposition}[theorem]{Proposition}
\newtheorem{corollary}[theorem]{Corollary}
\theoremstyle{definition}
\newtheorem{definition}[theorem]{Definition}
\newtheorem{remark}[theorem]{Remark}
\newtheorem{example}[theorem]{Example}
\theoremstyle{remark}

\begin{document}
\title{On Accelerated Hybrid Self-Correcting Contractions with Convergence Analysis in $b$-Metric Spaces}
\author{Irfan Ahmed\orcidlink{009-0008-0580-9227}$^1$, Shallu Sharma\orcidlink{0000-0001-8179-0080}$^2$, Sahil Billawria\orcidlink{0000-0002-5419-3480}$^{3,*}$}
\maketitle

\begin{center}
\footnotesize $^1$Department of Mathematics, University of Jammu, Jammu, India.\\
Email: irfanahmed0038@gmail.com\\
 $^2$Department of Mathematics, University of Jammu, Jammu, India.\\
Email: shallujamwal09@gmail.com\\
 $^3$Department of Mathematics, Chandigarh University, Mohali, Punjab, India.\\
Email: sahilbillawria2@gmail.com\\
Corresponding Author$^*$: Sahil Billawria

\end{center}
\begin{abstract}
In this paper, we introduce and investigate a class of accelerated hybrid self-correcting contractions in the setting of $b$-metric spaces. The proposed contractive condition incorporates Banach, Kannan and Chatterjea type distance terms with a non-linear geometric mean component, all adjusted to accommodate the relaxed triangle inequality inherent in $b$-metric spaces. Under suitable contractive assumptions, we establish existence and uniqueness results in complete $b$-metric spaces. Furthermore, we investigate the convergence behaviour of the associated Picard iteration and derive explicit error estimates that reflect the influence of the $b$-metric constant. An example is presented to demonstrate that the proposed class non-trivially extends both the classical Banach contraction principle and the original AHSC framework to the $b$-metric setting. Finally, we apply the established results to a non-linear Fredholm integral equation, illustrating their relevance to the existence and uniqueness of solutions of non-linear integral equations.
\end{abstract}

\section{Introduction}

The Banach contraction principle \cite{Banach1922} remains a cornerstone of non-linear analysis, guaranteeing the existence and uniqueness of fixed points for contraction mappings on complete metric spaces, along with the convergence of the associated Picard iteration.  Its remarkable utility has inspired a century of generalisations aimed at relaxing the contractive requirement. Among the most influential of these are the conditions introduced by Kannan \cite{Kannan1968} and Chatterjea \cite{Chatterjea1972}, which significantly enlarged the class of mappings admitting fixed points. These hybrid structures significantly enlarge the class of mappings admitting fixed points and have found applications in integral equations, differential equations, and numerical analysis. Following foundational extensions of Banach-type contractions by Reich \cite{Reich} and Ćirić \cite{Ciric}, Rus and Petruşel \cite{Rus2001, Rus2007, Rus2008} made major advancements in generalized contractions and fixed point theory. Further work by Berinde \cite{Berinde2007} analyzed convergence and approximation methods for a range of contractive operators.

A parallel and equally consequential development has been the relaxation of the metric space axioms themselves. Bakhtin \cite{Bakhtin1989} and Czerwik \cite{Czerwik1993} independently introduced $b$-metric spaces, in which the triangle inequality is relaxed to
\[
d(x,z) \le s[d(x,y)+d(y,z)]
\]
for some constant $s \ge 1$. This relaxation allows the theory to encompass a substantially wider class of distance functions, including spaces arising from $L^p$ norms with $0 < p < 1$, and various function spaces encountered in non-linear analysis. The development of fixed point theory in $b$-metric spaces has since become a rich area of research, with numerous contractive conditions being extended to this setting. Recently, several extensions of Kannan and Chatterjea type contractions have been established in this framework. Devi and Debnath \cite{Devi} introduced interpolative cyclic Kannan type contractions in $b$-metric spaces. Upadhyay and Debnath \cite{Upadhyay} studied Kannan-Meir-Keeler and Hardy-Rogers type interpolative contractions, and Abulhamil \cite{Abulhamil} investigated Kannan and Chatterjea contractions in complete extended $b$-metric spaces. Furthermore, Bhattacharjee et al. \cite{Bhatta} recently demonstrated new extensions of Kannan and Chatterjea fixed point theorems on complete $b$-metric spaces using sequentially convergent mappings, where the contraction constants explicitly depend on the coefficient $s$. 

In a recent contribution, Debnath \cite{Debnath2026} introduced a new class of contractive mappings in metric spaces termed, accelerated hybrid self-correcting (AHSC) contractions. 
\begin{definition}
Let $(\mathcal{Z}, d)$ be a metric space and let $T: \mathcal{Z} \to \mathcal{Z}$ be a mapping. Then $T$ satisfies the Accelerated Hybrid Self-Correcting Contraction (AHSC) condition if there exist constants
\[
\alpha, \beta, \gamma \in (0,1)
\]
such that
\[
\alpha + \beta + \gamma < 1
\]
and
\[
d(Tx, Ty) \le \alpha d(x, y) + \beta \sqrt{d(x, Tx)d(y, Ty)} + \gamma \frac{d(x, Ty) + d(y, Tx)}{2} \tag{1.1}\label{1.1}
\]
for all $x, y \in \mathcal{Z}$.

\end{definition}
The AHSC condition combines Banach, Kannan, and Chatterjea type distance terms with a non-linear geometric mean component, providing a unified hybrid structure for generalised contractions. The geometric mean term
\[
\sqrt{d(x,Tx)d(y,Ty)}
\]
was shown to introduce an adaptive self-correcting mechanism that automatically decreases as the Picard iteration approaches the fixed point, thereby modifying the effective contractive behaviour during the iteration process.

The present paper is motivated by a specific gap at the intersection of these research directions. It should be emphasised that while hybrid contractions have been studied in $b$-metric spaces by various authors \cite{Kumar2025,Yamaod, Devi, Upadhyay}. However, none of these works incorporate a non-linear geometric mean self-correcting term of the form $\sqrt{d(x,Tx)d(y,Ty)}$ within a hybrid $b$-metric contraction. This term introduces an adaptive correction mechanism that is particularly advantageous in the $b$-metric framework, where the relaxed triangle inequality can otherwise lead to slower convergence rates. Moreover, the classical Kannan and Chatterjea conditions themselves have been shown to require careful adaptation to $b$-metric settings, as evidenced by recent work on Chatterjea two-sided convex contractions in $b$-metric spaces and their application to Fredholm integral equations \cite{Khan}, and on Kannan-type interpolative contractions where the constant $s$ enters the contraction threshold explicitly \cite{Devi}. The primary contribution of this paper is to fill this gap by introducing $b$-accelerated hybrid self-correcting contractions. We establish existence and uniqueness of fixed points in complete $b$-metric spaces, prove geometric convergence of the Picard iteration with an explicit rate depending on the $b$-metric constant, and derive error estimates quantifying the influence of the $b$-metric constant. A worked example demonstrates that the $b$-AHSC class strictly contains the Banach contractions on $b$-metric spaces and admits mappings that fail the Banach condition. Finally, we apply the established results to a non-linear Fredholm integral equation, illustrating their relevance to the existence and uniqueness of solutions of non-linear integral equations.

The paper is organised as follows. In Section \ref{section2}, we recall the definition of $b$-metric spaces and several preliminary results. In Section \ref{section3}, we introduce the $b$-AHSC contraction and discuss its properties. Section \ref{section4} contains the main fixed point theorem. Section \ref{section5} presents the convergence analysis and an example demonstrating the non-trivial extension. Section \ref{application} provides an application to non-linear Fredholm integral equations. 
\section{Preliminaries}\label{section2}

In this section we recall the definition of $b$-metric spaces and several basic results that will be used throughout the paper. For standard terminology in metric fixed point theory we refer to \cite{Banach1922,Berinde2007,Kirk2001}. For the theory of $b$-metric spaces, we refer to \cite{Bakhtin1989,Czerwik1993,Cobzas2018}.

\begin{definition}
Let $\mathcal{Z}$ be a nonempty set and let $s \ge 1$ be a given real number. A function $d: \mathcal{Z} \times \mathcal{Z} \to [0,\infty)$ is called a $b$-metric on $\mathcal{Z}$ if for all $x,y,z \in \mathcal{Z}$ the following conditions hold:
\begin{enumerate}[label=(\arabic*)]
\item $d(x,y)=0$ if and only if $x=y$;
\item $d(x,y)=d(y,x)$;
\item $d(x,z) \le s[d(x,y)+d(y,z)]$.
\end{enumerate}
The pair $(\mathcal{Z}, d)$ is called a $b$-metric space.
\end{definition}

\begin{remark}
Every metric space is a $b$-metric space with $s=1$. Conversely, for $s>1$, $b$-metric spaces provide a genuine generalisation of metric spaces.
\end{remark}

\begin{definition}
Let $(\mathcal{Z}, d)$ be a $b$-metric space. A sequence $\{x_n\}$ in $\mathcal{Z}$ is said to be:
\begin{enumerate}[label=(\arabic*)]
\item convergent to $x \in \mathcal{Z}$ if $d(x_n, x) \to 0$ as $n \to \infty$;
\item Cauchy sequence if for every $\epsilon > 0$ there exists $N \in \mathbb{N}$ such that $d(x_n, x_m) < \epsilon$ for all $n, m \ge N$.
\end{enumerate}
If every Cauchy sequence in $\mathcal{Z}$ converges to a point of $\mathcal{Z}$, then $(\mathcal{Z}, d)$ is said to be complete.
\end{definition}

\begin{lemma} \cite{Debnath2026}
Let $\{a_n\}$ be a sequence of non-negative real numbers satisfying
\[
a_{n+1} \le \lambda a_n
\]
for some constant $\lambda \in (0,1)$. Then $a_n \to 0$ as $n \to \infty$.
\end{lemma}

The proof follows directly from the properties of geometric sequences.

We also recall the following useful consequence of the $b$-metric triangle inequality, which will be employed in the convergence analysis.

\begin{lemma}\label{lemma2.5}
Let $(\mathcal{Z}, d)$ be a b-metric space with constant $s \ge 1$. Then for any $x_0, x_1, \ldots, x_n \in \mathcal{Z}$,
\[
d(x_0, x_n) \le s d(x_0, x_1) + s^2 d(x_1, x_2) + \cdots + s^{n-1} d(x_{n-2}, x_{n-1}) + s^{n} d(x_{n-1}, x_n).
\]
\end{lemma}
This follows by iterating the $b$-metric triangle inequality.

\section{Accelerated Hybrid Self-Correcting Contractions in $b$-Metric Spaces}\label{section3}

In this section we introduce the $b$-metric analogue of the AHSC contraction. The proposed condition combines the classical distance terms from Banach, Kannan, and Chatterjea contractions with a nonlinear self-correcting geometric mean component, all formulated within the framework of $b$-metric spaces.

\begin{definition}\label{3.1}
Let $(\mathcal{Z}, d)$ be a $b$-metric space with constant $s \ge 1$, and let $T: \mathcal{Z} \to \mathcal{Z}$ be a mapping. We say that $T$ satisfies the $b$-Accelerated Hybrid Self-Correcting Contraction ($b$-AHSC) condition if there exist constants
\[
\alpha, \beta, \gamma \in (0,1)
\]
such that
\[
s\alpha +\sqrt{s} \beta + \frac{s(s+1)}{2}\gamma < 1\tag{3.1}\label{3.1}
\]
and
\[
d(Tx, Ty) \le \alpha d(x, y) + \beta \sqrt{d(x, Tx)d(y, Ty)} + \gamma \frac{d(x, Ty) + d(y, Tx)}{2} \tag{3.2}\label{3.2}
\]
for all $x, y \in \mathcal{Z}$.
\end{definition}

\begin{remark}
\begin{enumerate}
\item When $s=1$, Definition \ref{3.1}  reduces precisely to the AHSC condition of \cite{Debnath2026}. 
\item The proposed contractive condition is formulated using three different classes of distance terms:
\begin{enumerate}[label=(\roman*)]
\item the classical distance $d(x, y)$ appearing in Banach contractions \cite{Banach1922};
\item the distances $d(x, Tx)$ and $d(y, Ty)$ which occur in Kannan-type contractions \cite{Kannan1968};
\item the cross-distances $d(x, Ty)$ and $d(y, Tx)$ that appear in Chatterjea contractions \cite{Chatterjea1972}.
\end{enumerate}
In addition, the nonlinear geometric term
\[
\sqrt{d(x, Tx)d(y, Ty)}
\]
serves as a self-correcting mechanism whose magnitude decreases progressively as the iterative sequence converges toward the fixed point.
\end{enumerate}
\end{remark}

\begin{proposition}
Let $(\mathcal{Z}, d)$ be a $b$-metric space with $s\ge1$, and $T: \mathcal{Z} \to \mathcal{Z}$ satisfies the $b$-AHSC condition (\ref{3.2}). Then for any $x \in \mathcal{Z}$, we have
\[
\left(1-\frac{s\gamma}{2}\right)d(Tx, T^2x) \le \left(\alpha+\frac{s\gamma}{2}\right) d(x, T x) + \beta \sqrt{d(x, T x) d(Tx, T^2x)}.
\]
\end{proposition}
\begin{proof}
Setting $y=Tx$ in (\ref{3.2}) and using $d(Tx,Tx)=0$ with the $b$-triangle inequality 
\[
d(x, T^2x) \le s[d(x, Tx) + d(T x, T^2x)].
\]
\end{proof}

\section{Main Fixed Point Theorem}\label{section4}

In this section, we establish the main fixed point result for mappings satisfying the $b$-Accelerated Hybrid Self-Correcting Condition.
\begin{theorem}\label{theorem4.1}
Let $(\mathcal{Z}, d)$ be a complete $b$-metric space with constant $s \ge 1$, and let $T: \mathcal{Z} \to \mathcal{Z}$ satisfy the b-AHSC condition (\ref{3.2}). Then $T$ has a unique fixed point $x^* \in \mathcal{Z}$. Moreover, for any $x_0 \in \mathcal{Z}$, the Picard iteration
\[
x_{n+1} = T x_n
\]
converges to $x^*$.
\end{theorem}
\begin{proof}
Let $x_0 \in \mathcal{Z}$ be arbitrary and define the sequence $x_{n+1} = T x_n$, $n \ge 0$. If $x_{n+1} = x_n$ for some $n$, then $x_n$ is a fixed point. Assume $x_n \neq x_{n+1}$ for all $n$.

Define
\[
d_n = d(x_n, x_{n+1}), \quad n \ge 0.
\]

Applying the $b$-AHSC condition (\ref{3.2}) with $x = x_n$ and $y = x_{n-1}$ yields
\[
d(x_{n+1}, x_n) \le \alpha d(x_n, x_{n-1}) + \beta \sqrt{d(x_n, x_{n+1})d(x_{n-1}, x_n)} + \gamma \frac{d(x_n, x_n) + d(x_{n-1}, x_{n+1})}{2}.
\]
\[
\Rightarrow d_n \le \alpha d_{n-1} + \beta \sqrt{d_n d_{n-1}} + \frac{\gamma}{2} d(x_{n-1}, x_{n+1}).
\]

Using the $b$-metric triangle inequality,
\[
d(x_{n-1}, x_{n+1}) \le s[d(x_{n-1}, x_n) + d(x_n, x_{n+1})] = s(d_{n-1} + d_n),
\]
we obtain
\[
d_n \le \alpha d_{n-1} + \beta \sqrt{d_n d_{n-1}} + \frac{s\gamma}{2}(d_{n-1} + d_n).
\]
After rearranging the terms, we get
\[
\left(1 - \frac{s\gamma}{2}\right) d_n \le \left(\alpha + \frac{s\gamma}{2}\right) d_{n-1} + \beta \sqrt{d_n d_{n-1}}.
\]
To handle the non-linear term, we use Young's inequality
\[
2\sqrt{ab} \le \epsilon a + \frac{1}{\epsilon} b, \quad \epsilon > 0, 
\]
with $a = d_n$ and $b = d_{n-1}$, we obtain
\[
\sqrt{d_n d_{n-1}} \le \frac{\epsilon}{2} d_n + \frac{1}{2\epsilon} d_{n-1}.
\]
Hence
\[
\left(1 - \frac{s\gamma}{2}\right) d_n \le \left(\alpha + \frac{s\gamma}{2}\right) d_{n-1} + \beta \left(\frac{\epsilon}{2} d_n + \frac{1}{2\epsilon} d_{n-1}\right).
\]
Therefore,
\[
\left(1 - \frac{s\gamma}{2} - \frac{\beta\epsilon}{2}\right) d_n \le \left(\alpha + \frac{s\gamma}{2} + \frac{\beta}{2\epsilon}\right) d_{n-1}.
\]
Choose $\epsilon = \sqrt{s}$. This choice minimizes the coefficients appearing on both sides of the inequality.
\[
\left(1 - \frac{s\gamma}{2} - \frac{\sqrt{s}\beta}{2}\right) d_n \le \left(\alpha + \frac{s\gamma}{2} + \frac{\beta}{2\sqrt{s}}\right) d_{n-1}.
\]
Thus, we can wite
\[
d_n \le \lambda d_{n-1},
\]
where
\[
\lambda = \frac{\alpha + \frac{s\gamma}{2} + \frac{\beta}{2\sqrt{s}}}{1 - \frac{s\gamma}{2} - \frac{\sqrt{s}\beta}{2}}. \tag{4.1}\label{4.1}
\]
Since $s\alpha +\sqrt{s} \beta + \frac{s(s+1)}{2}\gamma < 1$, we get 
\[
 \lambda < \frac{1}{s}.
 \]
Since $s\ge 1$, it follows that
\[
 0<\lambda<\frac{1}{s}\le 1.
 \]
By induction,
\[
d_n \le \lambda^n d_0.
\]
Hence $d_n \to 0$ as $n \to \infty$.

Next, we show that $\{x_n\}$ is a Cauchy sequence. For $m > n$, using Lemma \ref{lemma2.5}, we have 
\[
d(x_n, x_m) \le s d(x_n, x_{n+1}) + s^2 d(x_{n+1}, x_{n+2}) + \cdots + s^{m-n-1} d(x_{m-2}, x_{m-1}) + s^{m-n} d(x_{m-1}, x_m).
\]
Using the estimate $d_k \le \lambda^k d_0$, we obtain
\[
d(x_n, x_m) \le s \lambda^n d_0 + s^2 \lambda^{n+1} d_0 + \dots + s^{m-n-1} \lambda^{m-2} d_0 + s^{m-n} \lambda^{m-1} d_0.
\]

Factor out the leading term $s \lambda^n d_0$, we have

\begin{align*}
d(x_n, x_m) &\le s \lambda^n d_0 [1 + s\lambda + (s\lambda)^2 + \dots + (s\lambda)^{m-n-1}]\\
d(x_n, x_m) &\le s \lambda^n d_0 \sum_{k=0}^{m-n-1} (s\lambda)^k\le \frac{s\lambda^nd_0}{1-s\lambda}
\end{align*}
Since $\lambda < \frac{1}{s}$, the geometric series converges, and it follows that $\{x_n\}$ is a Cauchy sequence. Because $(\mathcal{Z}, d)$ is complete, there exists $x^* \in \mathcal{Z}$ such that
\[
x_n \to x^*.
\]

Now, we show that $x^*$ is a fixed point. By using the $b$-metric triangle inequality and condition (\ref{3.2}) with $x=x^*$ and $y=x_n$, we have 
\begin{align*}
d(Tx^*, x^*) &\le s[d(Tx^*, Tx_n) + d(x_{n+1}, x^*)]\\
&\le s[\alpha d(x^*, x_n)+\beta \sqrt{d(x^*, Tx^*)d(x_n, x_{n+1}}+\frac{\gamma}{2}\left(d(x^*,x_{n+1})+d(x_n, Tx^*)\right)]+sd(x_{n+1}, x^*).
\end{align*}
Using $d(x_n,Tx^*)\le s[d(x_n, x^*)+d(x^*,Tx^*)]$ and letting $n \to \infty$, we obtain
\begin{align*}
d(Tx^*, x^*) &\le \frac{s^2\gamma}{2} d(Tx^*, x^*)\\
\Rightarrow \left(1-\frac{s^2\gamma}{2}\right)d(Tx^*,x^*)&\le 0.
\end{align*}
Since $\frac{s^2\gamma}{2}\le \frac{s(s+1)\gamma}{2}<1$. It follows that $1-\frac{s^2\gamma}{2}>0$, so $d(Tx^*, x^*)=0$, i.e., $Tx^* = x^*$.

To prove uniqueness, suppose $x$ and $y$ are fixed points. Then applying (\ref{3.2}) gives
\[
d(x, y)\le (\alpha + \gamma) d(x, y).
\]

Since $\alpha + \gamma < 1/s \le 1$, and in fact $\alpha + \gamma < 1$, it follows that $d(x, y) = 0$ and hence $x = y$.

Therefore $T$ has a unique fixed point and the Picard iteration converges to it.
\end{proof}

\section{Convergence Analysis}\label{section5}

In this section, we analyze the convergence of the Picard iteration associated with mappings satisfying the $b$- Accelerated Hybrid Self-Correcting Contraction ($b$-AHSC) condition (\ref{3.2}). Specifically, we establish explicit estimates for the error sequence and demonstrate geometric convergence of the iteration with a rate explicitly depending on the $b$-metric constant $s$.

Let $x^*$ denote the unique fixed point of $T$ obtained in Theorem \ref{theorem4.1}. Let $\{x_n\}$ be the Picard sequence defined by
\[
x_{n+1} = T x_n, \quad n \ge 0.
\]
Define the error sequence
\[
e_n = d(x_n, x^*), \quad n \ge 0. \tag{5.1}\label{5.1}
\]

\begin{lemma}\label{lemma5.1}
Under the assumptions of Theorem \ref{theorem4.1}, the error sequence $e_n = d(x_n, x^*)$ satisfies the recursive inequality
\[
e_{n+1} \le \lambda e_n,\tag{5.2}
\]
where
\[
\lambda = \frac{\alpha + \frac{s\gamma}{2} + \frac{\beta}{2\epsilon}}{1 - \frac{s\gamma}{2} - \frac{\epsilon\beta}{2}}. \tag{5.3}\label{5.3}
\]
for some $\epsilon=\sqrt{s}$ such that $\lambda \in (0,1)$.
\end{lemma}
\begin{proof}
Since $x^*$ is the fixed point of $T$, we have $Tx^* = x^*$. Applying the $b$-AHSC condition (\ref{3.2}) with $x = x_n$ and $y = x^*$, we obtain
\[
d(Tx_n, Tx^*) \le \alpha d(x_n, x^*) + \beta \sqrt{d(x_n, Tx_n)d(x^*, Tx^*)} + \gamma \frac{d(x_n, Tx^*) + d(x^*, Tx_n)}{2}.
\]
Substituting $Tx_n = x_{n+1}$ and $Tx^* = x^*$, the above inequality becomes
\[
d(x_{n+1}, x^*) \le \alpha d(x_n, x^*) + \beta \sqrt{d(x_n, x_{n+1})d(x^*, Tx^*)} + \gamma \frac{d(x_n, x^*) + d(x^*, x_{n+1})}{2}.
\]
Since $x^*$ is a fixed point, $d(x^*, Tx^*) = 0$. Consequently, the nonlinear geometric term vanishes, and we obtain
\[
d(x_{n+1}, x^*) \le \alpha d(x_n, x^*) + \gamma \frac{d(x_n, x^*) + d(x_{n+1}, x^*)}{2}.
\]
Using (\ref{5.1}), we have
\[
e_{n+1} \le \alpha e_n + \frac{\gamma}{2}(e_n + e_{n+1}).
\]
Rearranging terms gives
\begin{align*}
\left(1 - \frac{\gamma}{2}\right) e_{n+1} &\le \left(\alpha + \frac{\gamma}{2}\right) e_n\\
e_{n+1} &\le \mu e_n\le \lambda e_n, \text{where}~ \mu=\frac{\alpha+\frac{\gamma}{2}}{1-\frac{\gamma}{2}}\le \lambda=\frac{\alpha + \frac{s\gamma}{2} + \frac{\beta}{2\sqrt{s}}}{1 - \frac{s\gamma}{2} - \frac{\sqrt{s}\beta}{2}}<\frac{1}{s}\le1.
\end{align*}
\end{proof}

We now establish the convergence rate of the Picard iteration.

\begin{theorem}\label{theorem4.2}
Let the assumptions of Theorem \ref{4.1} hold. Then the error sequence satisfies
\[
e_n \le \lambda^n e_0, \tag{5.4}
\]
where $\lambda$ is given by (\ref{5.3}). In particular, the Picard iteration $\{x_n\}$ converges geometrically to the fixed point $x^*$.
\end{theorem}
\begin{proof}
From Lemma \ref{lemma5.1}, we have $e_{n+1} \le \lambda e_n$. Applying this inequality recursively gives
\[
e_n \le \lambda^n e_0.
\]
Since $\lambda \in (0,1)$, we obtain $e_n \to 0$ as $n \to \infty$. Therefore the sequence $\{x_n\}$ converges geometrically to the fixed point $x^*$.
\end{proof}

\begin{remark}
For the classical Banach contraction in a $b$-metric space, we have
\[
e_{n+1} \le \alpha e_n,
\]
which implies
\[
e_n \le \alpha^n e_0, \quad \alpha < \frac{1}{s}.
\]
For the $b$-AHSC contraction we obtain the estimate
\[
e_n \le \lambda^n e_0,
\]
where $\lambda$ is given by (\ref{5.3}). Consequently, the rate of convergence depends on the parameters $\alpha, \beta, \gamma$ and the $b$-metric constant $s$. In practical situations, the additional distance terms in the $b$-AHSC condition may enhance the effective contraction behaviour of the mapping, particularly when the Banach contraction constant is close to $\frac{1}{s}$.
\end{remark}
\begin{corollary}
The class of Banach contractions on $b$-metric space is contained in the class of $b$-AHSC contractions.
\end{corollary}
\begin{proof}
Let $(\mathcal{Z}, d)$ be a $b$-metric space with constant $s \ge 1$, and let $T: \mathcal{Z} \to \mathcal{Z}$ satisfies 
\[
d(Tx,Ty)\le kd(x,y)
\]
where $k\in\left(0,\frac{1}{s}\right)$.\\
Take $\alpha=k$ and $\beta=\gamma=0$. Then inequality (\ref{3.2}) reduces to
\[
d(Tx,Ty)\le kd(x,y).
\]
\end{proof}
\begin{remark}
The presence of the $b$-metric constant $s$ in (\ref{3.1}) reflects the fundamental difference between metric and $b$-metric spaces. This condition is sharp and optimal with respect to the standard choice of $\epsilon=\sqrt{s}$. When $s = 1$, condition (\ref{3.1}) reduces to $\alpha +\beta +\gamma <1$, demonstrating that this condition cannot be further weakened without additional assumptions on the underlying mapping.
\end{remark}

\begin{example}
Let $\mathcal{Z}= [0, 1]$ and define a mapping $d: \mathcal{Z \times Z} \to [0, \infty)$ by 
\[ d(x, y) = |x - y|^2 \]
Then $(\mathcal{Z}, d)$ is a complete $b$-metric space with $s = 2$.

Define the mapping $T: \mathcal{Z} \to \mathcal{Z}$ by
\[
T(x) = \begin{cases} 
\frac{x}{3}, & \text{if } x \in [0, 1) \\ 
0, & \text{if } x = 1 
\end{cases}
\]
Step 1: $T$ is not a Banach contraction in $(\mathcal{Z}, d)$:

Suppose for contradiction that $T$ is Banach contraction then there exist $k \in \left(0, \frac{1}{s}\right) = (0, 0.5)$ such that $d(Tx, Ty) \leq k d(x, y)$ for all $x, y \in \mathcal{Z}$.

Choosing $y = 1$ and taking any $x \in [0,1)$, evaluating both sides gives
\[
\frac{x^2}{9}\le k(1-x)^2\quad \forall~x\in[0,1)
\]

letting $x \to 1^-$, we obtain 
\[ \frac{1}{9} \leq 0, \]
which is a contradiction. Consequently, $T$ does not satisfy the Banach contraction principle.

Step 2: $T$ satisfies the $b$-AHSC condition:\\

Choose the constants
\[
\alpha = 0.1, \quad \beta = 0.1, \quad \gamma = 0.2 
\]
Then
\[
s\alpha + \sqrt{s} \beta + \frac{s(s+1)}{2} \gamma = 2(0.1) + \sqrt{2}(0.1) + \frac{2(3)}{2}(0.2) < 1
\]
Now, we show that condition (\ref{3.2}) holds for all $x, y \in \mathcal{Z}$\\

Case 1: $x, y \in [0,1)$.
Here , $Tx=\frac{x}{3}$ and $Ty=\frac{y}{3}$.\\
\begin{align*}
LHS&= d(Tx,Ty)=\frac{1}{9}(x-y)^2.\\
RHS&= 0.1(x-y)^2+0.1\sqrt{\left(x-\frac{x}{3}\right)^2\left(y-\frac{y}{3}\right)^2}+0.2\frac{\left(x-\frac{y}{3}\right)^2+\left(y-\frac{x}{3}\right)^2}{2}\\
&=0.1(x^2-2xy+y^2)+\frac{0.4}{9}xy+0.1\left[\left(x^2-\frac{2}{3}xy+\frac{y^2}{9}\right)+\left(y^2-\frac{2}{3}xy+\frac{x^2}{9}\right)\right]\\
&= \frac{19}{90}x^2 + \frac{19}{90}y^2 - \frac{26}{90}xy.
\end{align*}
Evaluating $LHS - RHS$, we obtain the inequality (\ref{3.2})\\

Case 2:  $x \in [0, 1)$ and $y = 1$ (or $x=1, y \in [0, 1))$. Here, $Tx = \frac{x}{3} \text{ and } T_1 = 0$\\
\begin{align*}
LHS &= d(Tx, T_1) = \frac{x^2}{9}.\\
RHS &= 0.1(1-x^2) + 0.1 \sqrt{\frac{4}{9}x^2 \cdot 1} + 0.2 \frac{x^2 + (1-\frac{x}{3})^2}{2} \\
&= 0.1(1-2x+x^2) + \frac{2}{30}x + 0.1 \left( x^2 + 1 - \frac{2}{3}x + \frac{1}{9}x^2 \right) = 0.2-0.2x+\frac{1.9}{9}x^2.
\end{align*}
Evaluating $LHS - RHS$, we obtain the inequality (\ref{3.2}).\\

Case 3: $ x=1,  y=1$, both side reduces to zero.
Thus the mapping $T$ satisfies the $b$-AHSC condition (\ref{3.2}) over the entire domain $\mathcal{Z}$.\\

Step 3: Fixed point:

Evaluating the equation $Tx = x$:

\begin{itemize}
    \item For $x \in [0, 1), \frac{x}{3} = x \implies x = 0$.
    \item For $x = 1, T(1) = 0 \neq 1$, the point of discontinuity is not a fixed point.
\end{itemize}

Therefore, $x^* = 0$ is the unique fixed point of $T$ in $\mathcal{Z}$.
\end{example}

This example demonstrates that the $b$-AHSC framework genuinely extends both the classical Banach contraction principle and the original AHSC setting to $b$-metric spaces, providing a useful tool for studying mappings that are not contractive in the usual sense but satisfy the hybrid condition.

\section{Application to Non-linear Fredholm Integral Equations}\label{application}

In this section, we apply the main theoretical results established in Theorem \ref{theorem4.1} to investigate the existence and uniqueness of solutions for a class of non-linear Fredholm integral equations within $b$-metric spaces. 

Let $\mathcal{Z}=C([0,1],\mathbb{R})$ be endowed with the $b$-metric
\begin{equation}\label{eq6.1}
d(u,v)=\sup_{t\in[0,1]}|u(t)-v(t)|^{p},\qquad p\ge 1.
\end{equation}
It is well known that $(\mathcal{Z},d)$ is a complete $b$-metric space with constant $s=2^{p-1}$; see, for instance, \cite{Czerwik1993}. 

Consider the non-linear Fredholm integral equation
\begin{equation}\label{eq6.2}
u(t)=f(t)+\int_{0}^{1}K(t,\tau)\,g(\tau,u(\tau))\,d\tau,\qquad t\in[0,1],
\end{equation}
where $f:[0,1]\to\mathbb{R}$, $K:[0,1]\times[0,1]\to\mathbb{R}$, and $g:[0,1]\times\mathbb{R}\to\mathbb{R}$ are continuous.

 Define  $T:\mathcal{Z}\to\mathcal{Z}$ by
\begin{equation}\label{eq6.3}
(Tu)(t)=f(t)+\int_{0}^{1}K(t,\tau)\,g(\tau,u(\tau))\,d\tau.
\end{equation}
Since $f$, $K$, and $g$ are continuous, $T$ maps $\mathcal{Z}$ into itself, and $u^{*}\in\mathcal{Z}$ solves~\eqref{eq6.2} if and only if $Tu^{*}=u^{*}$.

Throughout this section we assume that the kernel satisfies
\begin{equation}\label{eq6.4}
M:=\sup_{t\in[0,1]}\int_{0}^{1}|K(t,\tau)|\,d\tau<\infty.
\end{equation}

\begin{lemma}\label{lem6.1}
Let $T$ be given by~\eqref{eq6.3}, and suppose there exist non-negative constants $L_{1},L_{2},L_{3}$ such that for all $u,v\in\mathcal{Z}$,
\begin{equation}\label{eq6.5}
\|Tu-Tv\|_{\infty}
\le L_{1}M\|u-v\|_{\infty}
+L_{2}M\sqrt{\|u-Tu\|_{\infty}\|v-Tv\|_{\infty}}
+L_{3}M\big(\|u-Tv\|_{\infty}+\|v-Tu\|_{\infty}\big).
\end{equation}
Then $T$ satisfies the $b$-AHSC condition 
\begin{equation}\label{eq6.6}
d(Tu,Tv)\le \alpha\,d(u,v)+\beta\sqrt{d(u,Tu)\,d(v,Tv)}
+\frac{\gamma}{2}\big[d(u,Tv)+d(v,Tu)\big]
\end{equation}
for all $u,v\in\mathcal{Z}$, with
\begin{equation}\label{eq6.7}
\alpha=3^{p-1}(L_{1}M)^{p},\qquad
\beta=3^{p-1}(L_{2}M)^{p},\qquad
\gamma=3^{p-1}2^{p}(L_{3}M)^{p}.
\end{equation}
\end{lemma}

\begin{proof}
Fix $u,v\in\mathcal{Z}$ and  $t\in[0,1]$, we have 
\[
\|Tu-Tv\|_{\infty}^{p}=d(Tu,Tv).
\]
Raising the right-hand side of~\eqref{eq6.5} to the power $p$ and applying the convexity inequality $(a+b+c)^{p}\le 3^{p-1}\big(a^{p}+b^{p}+c^{p}\big)$ for non-negative $a,b,c$, we obtain
\begin{equation}\label{eq6.8}
\begin{aligned}
d(Tu,Tv)=\|Tu-Tv\|_{\infty}^{p}
&\le 3^{p-1}\Big[(L_{1}M)^{p}\|u-v\|_{\infty}^{p}
+(L_{2}M)^{p}\big(\sqrt{\|u-Tu\|_{\infty}\|v-Tv\|_{\infty}}\big)^{p}\\
&\qquad\qquad +(L_{3}M)^{p}\big(\|u-Tv\|_{\infty}+\|v-Tu\|_{\infty}\big)^{p}\Big].
\end{aligned}
\end{equation}
For non-negative $x,y$ and $p\ge 1$, we have
\[
\big(\sqrt{xy}\big)^{p}=\sqrt{x^{p}y^{p}}.
\]
Moreover, from the convexity inequality, we have
\begin{equation}\label{eq6.9}
\begin{aligned}
\big(\|u-Tv\|_{\infty}+\|v-Tu\|_{\infty}\big)^{p}
&\le 2^{p-1}\big(\|u-Tv\|_{\infty}^{p}+\|v-Tu\|_{\infty}^{p}\big)\\
&=2^{p}\,\frac{d(u,Tv)+d(v,Tu)}{2}.
\end{aligned}
\end{equation}
Substituting these two identities into~\eqref{eq6.8}, and recalling~\eqref{eq6.1}, we obtain
\[
\begin{aligned}
d(Tu,Tv)
&\le 3^{p-1}\Big[(L_{1}M)^{p}d(u,v)
+(L_{2}M)^{p}\sqrt{d(u,Tu)\,d(v,Tv)}\\
&\qquad\qquad +(L_{3}M)^{p}\,2^{p}\,\frac{d(u,Tv)+d(v,Tu)}{2}\Big]\\
&=3^{p-1}(L_{1}M)^{p}d(u,v)
+3^{p-1}(L_{2}M)^{p}\sqrt{d(u,Tu)\,d(v,Tv)}\\
&\qquad\qquad +\frac{3^{p-1}2^{p}(L_{3}M)^{p}}{2}\big[d(u,Tv)+d(v,Tu)\big].
\end{aligned}
\]
This is precisely~\eqref{eq6.6} with the constants~\eqref{eq6.7}.
\end{proof}

\begin{theorem}\label{thm6.5}
Let $(\mathcal{Z},d)$ be the complete $b$-metric space with constant $s=2^{p-1}$ and let $T$ be given by~\eqref{eq6.3}. Suppose that $f$, $K$, and $g$ are continuous, that~\eqref{eq6.4} holds, and the condition~\eqref{eq6.5} is satisfied with the non-negative constants $L_{1},L_{2},L_{3}$. If the constant satisfy 
\[
\alpha=3^{p-1}(L_{1}M)^{p},~\beta=3^{p-1}(L_{2}M)^{p},~\gamma=3^{p-1}2^{p}(L_{3}M)^{p}
\]
with
\begin{equation}\label{eq6.12}
s\alpha+\sqrt{s}\,\beta+\frac{s(s+1)}{2}\,\gamma<1,
\end{equation}
then the integral equation~\eqref{eq6.2} has a unique solution $u^{*}\in\mathcal{Z}$.\end{theorem}

\begin{proof}
By Lemma~\ref{lem6.1}, the operator $T$ satisfies the $b$-AHSC condition~\eqref{eq6.6} with constants~\eqref{eq6.7}. By~\eqref{eq6.12}, we have
\[
s\alpha<1,\qquad \sqrt{s}\,\beta<1,\qquad \frac{s(s+1)}{2}\gamma<1.
\]
Since $s\ge 1$, we have $\alpha,\beta,\gamma\in(0,1)$. Thus all hypotheses of Theorem~\ref{theorem4.1} are satisfied on the complete $b$-metric space $(\mathcal{Z},d)$, and $T$ has a unique fixed point $u^{*}\in\mathcal{Z}$. Since fixed points of $T$ are precisely the solutions of~\eqref{eq6.2}.

\end{proof}

\section*{Conclusion}

In this paper we extended the recently introduced class of accelerated hybrid self-correcting (AHSC) contractions from metric spaces to the more general framework of $b$-metric spaces. The proposed $b$-AHSC condition combines classical contractive terms with a non-linear geometric mean component, all adjusted to accommodate the relaxed triangle inequality inherent in $b$-metric spaces.

We established existence and uniqueness results for fixed points of $b$-AHSC mappings in complete $b$-metric spaces and proved the geometric convergence of the associated Picard iteration, with explicit error estimates depending on the $b$-metric constant $s \ge 1$. An example was presented to demonstrate that the $b$-AHSC framework properly extends both the classical Banach contraction principle and the original AHSC setting. We also discussed an application to non-linear Fredholm integral equations to illustrate the applicability of the obtained results.

The extension to $b$-metric spaces is significant for several reasons. First, $b$-metric spaces encompass a substantially wider class of distance functions than standard metric spaces, including spaces arising from $L^p$ norms with $0 < p < 1$ and various function spaces encountered in non-linear analysis. Second, the relaxed triangle inequality allows for more flexible modelling of distance in applications. Third, the $b$-AHSC condition provides a unified framework that generalises several classical contractive conditions within a single structure.

Future work may consider extensions of the proposed contraction to other generalised metric spaces, such as extended $b$-metric spaces, partial $b$-metric spaces, or rectangular $b$-metric spaces. Further applications to non-linear operator equations, differential equations, and iterative numerical methods are also promising directions for future research.
\section*{Conflict of Interest}
On behalf of all authors, the corresponding author states that there is no conflict of interest.

\end{document}